%
\input amstex
\documentstyle{amsppt}
\magnification 1200
\NoBlackBoxes
\input epsf
\vcorrection{-12mm}

\def\Z{{\Bbb Z}}

\def\R{{\Bbb R}}
\def\eps{\varepsilon}
\def\ssim{\overset s\to\sim}
\def\End{\operatorname{End}}
\def\red{\operatorname{red}}
\def\sh{\operatorname{sh}}
\def\pr{\operatorname{pr}}
\def\Fred{\hat F^{\red}}
\def\rr{{\bold r}}

\def\sectHatB  {2}
\def\sectHe    {3}
\def\sectI     {4}

\def\thHatB     {\sectHatB.1}
\def\lemHatBone {\sectHatB.2}
\def\lemHatBtwo {\sectHatB.3}
\def\remFT      {\sectHatB.4}
\def\thMM       {\sectHe.1}
\def\thI        {\sectI.1}
\def\remBetaZero{\sectI.2}
\def\remGB      {\sectI.3}

\def\eqR        {1}
\def\eqCubRel   {2}
\def\eqDiagram  {3}

\def\figHatB {1}
\def\figYZ   {2}

\def\refBennequin {1}
\def\refBM        {2}
\def\refFT        {3}
\def\refFunar     {4}
\def\refJones     {5}
\def\refLNgS      {6}
\def\refNgOT      {7}
\def\refOrFu      {8}
\def\refOrGB      {9}
\def\refOSh       {10}
\def\refOST       {11}

\topmatter
\title
       Cubic Hecke Algebras and Invariants of Transversal Links
\endtitle
\author
       S.~Yu.~Orevkov
\endauthor
\address   Steklov Math. Institute, Moscow
\endaddress
\address   IMT, Universit\'e Paul Sabatier (Toulouse-3)
\endaddress
\email     orevkov\@math.ups-tlse.fr
\endemail

\endtopmatter

\document

\subhead\S1. Introduction
\endsubhead
Let $\alpha$ be a differential 1-form which defines the standard (tight) contact structure
in $\R^3$, e.~g., $\alpha=x\,dy-y\,dx+dz$.
A link $L$ in $\R^3$ is called {\it transversal} if $\alpha|_L$
does not vanish on $L$. Transversal links are considered up to isotopies such that
the link remains transversal at every moment. Transversal links and their invariants are being
actively studied, see, e.~g., [\refBM, \refLNgS, \refNgOT,
\refOST] and numerous references therein.
In the present paper, we propose a purely algebraic approach to construct
invariants of transversal links (similar to Jones' approach [\refJones] to construct
invariants of usual links). The only geometry used is the transversal analogue of
Alexander's and Markov's theorems proved in [\refBennequin] and [\refOSh] respectively.


Let $B_n$ be the group of $n$-braids.
We denote its standard (Artin's) generators by $\sigma_1,\dots,\sigma_{n-1}$.
Let $B_\infty=\lim B_n$ be the limit under the embeddings $B_n\to B_{n+1}$,
$\sigma_i\mapsto\sigma_i$.
Let $k$ be a commutative ring and $u,v$ indeterminates. We set $A=k[u]$,
$A_v=k[u,v]$ and we denote the corresponding group algebras by
$kB_\infty$, $AB_\infty$ and $A_vB_\infty$.
Let $\pi:kB_\infty\to H_\infty$ be a surjective morphism of $k$-algebras.
We extend it to the morphisms (also denotes by $\pi$) of $A$- and $A_v$-algebras
$AB_\infty\to AH_\infty=H_\infty\otimes_k A$ and
$A_vB_\infty\to A_vH_\infty=H_\infty\otimes_k A_v$.

Let $R$ be the $A$-submodule of $AB_\infty$
generated by all the elements of the form
$$
   XY-YX,\;\; X\sigma_n-uX\quad\text{ where }\; X,Y\in B_n,\; n\ge1,  \eqno(\eqR)
$$
and let $R_v$ be the $A_v$-submodule of $A_vB_\infty$ generated by (\eqR)
and also by $X\sigma_n^{-1}-vX$ for $X\in B_n$, $n\ge1$.
Let $M=AH_\infty/\pi(R)$ and $M_v = A_vH_\infty/\pi(R_v)$.
We say that the quotient map $t_v:A_vH_\infty\to M_v$ is the
{\it universal Markov trace} on $H_\infty$. Due to Alexander's and Markov's theorems,
it defines a link invariant
$P_{t_v}(L)=u^{(-n-e)/2}v^{(-n+e)/2}t_v(X)\in M_v\otimes_{A_v} k[u^{\pm1/2},v^{\pm1/2}]$
where $L$ is the closure of an $n$-braid $X$ and 
$e=e(X)=\sum_j e_j$ for $X=\prod_j \sigma_{i_j}^{e_j}$.

Similarly, by the transversal analogue of Alexander's and
Markov's theorems, the quotient map%
\footnote{I propose to call it {\it universal semi-Markov trace} on $H_\infty$.}
$t:AH_\infty\to M$ defines a transversal link invariant
$P_{t}(L)=u^{-n}t(X)\in M\otimes_A k[u^{\pm1}]$
where $L$ is the closure of an $n$-braid $X$.

Of course, these invariants do not make much sense unless there is a
reasonable solution to the identity problem in $M$ or in $M_v$. For example, if
$\ker\pi=0$, then $P_t$ is not really better than the tautological invariant
$I(L)=L$. However, if
$k=\Z[\alpha]$ and $A_vH_\infty=A_vB_\infty/(\sigma_1^2+\alpha\sigma_1+1)$, then
$M_v=A_v/(u+\alpha+v)\cong A$ and $P_{t_v}$ is the HOMFLY-PT polynomial up to
variable change.

In [\refOrFu], a description of $M_v$ is given when $H_\infty$
is a quotient of $kB_\infty$ by cubic relations of the form
$$
  \sigma_1^3-\alpha\sigma_1^2+\beta\sigma_1=1,\quad
  \sigma_2^\delta\sigma_1^{-\delta}\sigma_2^\delta=
  \sum_{\eps\in E} c_{\eps,\delta} \sigma_1^{\eps_1}\sigma_2^{\eps_2}\sigma_1^{\eps_3},
  \quad \delta=\pm1,
                                                                  \eqno(\eqCubRel)
$$
$E=\{\eps\in\{-1,0,1\}^3\mid\eps_2=0\Rightarrow\eps_1\eps_3=0\}$,
$\alpha,\beta,c_{\eps,\delta}\in k$. We have in this case $M_v=A_v/I_v$ and
a Gr\"obner base of $I_v$ can be computed at least theoretically.
Moreover, $I_v$ is computed in practice in a particular case when $H_\infty$
is the Funar algebra [\refFunar] and $\beta=0$. The computations may be fastened
using [\refOrGB].

In this paper we adapt the construction from [\refOrFu] for the computation of
$M$ when $H_\infty$ is defined by (\eqCubRel). We show that in this case $M\cong\hat A/\hat I$
where $\hat A=A[v_1,v_2,\dots\,]$ and $\hat I$ is an ideal of $\hat A$.
For any $d$, we give an algorithm to compute the ideal
$\hat I+(v_{d+1},v_{d+2},\dots)$.
Thus, we define an infinite sequence (indexed by $d$) of computable
transversal link invariants which carries the same information as the
universal semi-Markov trace on the cubic Hecke algebra given by (\eqCubRel).



\subhead\S\sectHatB. Monoid of braids with marked points
\endsubhead
Let $\hat B_n$ be the monoid of $n$-braids with a finite number of
points marked on the strings.
Algebraically it can be described as the monoid generated by
$\sigma_1^{\pm1},\dots,\sigma_{n-1}^{\pm1}$,$q_1,\dots,q_n$
(see Figure \figHatB),
subject to the standard braid group relations and the relations
$q_iq_j=q_jq_i$, $i,j=1,\dots,n-1$, and
$q_i\sigma_j = \sigma_j q_{T_j(i)}$ where $T_j$ is the transposition $(j,j+1)$.
Each element of $\hat B_n$ can be written in a unique way in the form
$q_1^{a_1}\dots q_n^{a_n}X$, $X\in B_n$, $a_i\ge 0$,
so, $\hat B_n=Q_n\rtimes B_n$ where $Q_n$ is
the free abelian monoid generated by $q_1,\dots,q_n$.

\midinsert
\epsfxsize=55mm
\centerline{\epsfbox{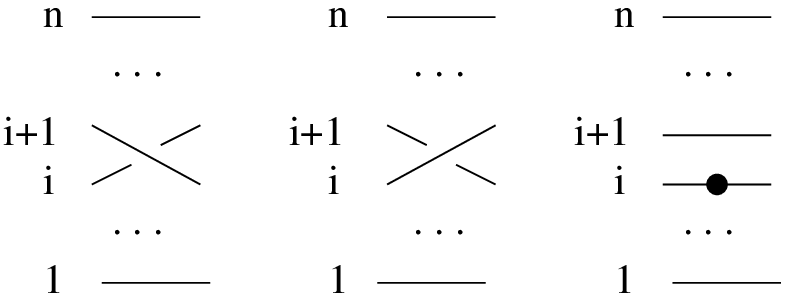}}
\centerline{\hskip 5mm $\sigma_i$ \hskip 15mm $\sigma_i^{-1}$ \hskip 13mm $q_i$}
\botcaption{Figure \figHatB} Generators of $\hat B_n$
\endcaption
\endinsert

Let $\hat B_\infty^\sqcup$ be the disjoint union $\bigsqcup_{n=1}^\infty\hat B_n$.
If an ambiguity is possible, we use the notation $(X)_n$ to emphasize that a word
$X$ represents an element of $\hat B_n$,
%
for example, the braid closure of $(1)_n$ is the trivial $n$-component link.

\proclaim{ Theorem \thHatB }
Transversal links are in bijection with the quotient of $\hat B_\infty^\sqcup$
by the equivalence relation generated by
$$
\xalignat 3
&(XY)_n\sim(YX)_n,  && X,Y\in\hat B_n,\; n\ge1 &&\text{{\rm(}conjugations{\rm)}},\\
&(X)_n \sim (X\sigma_n)_{n+1}, && X\in\hat B_n,\; n\ge1 &&
                    \text{{\rm(}positive Markov moves{\rm)}},\\
&(X q_n)_n \sim(X\sigma_n^{-1})_{n+1}, &&  X\in\hat B_n,\; n\ge1 &&
                    \text{{\rm(}negative Markov $q$-moves{\rm)}}\\
\endxalignat
$$
\endproclaim

\demo{ Proof } Follows easily from Lemma \lemHatBtwo.\qed\enddemo

Let $\ssim$ (strong equivalence) be the equivalence relation
on $\hat B_\infty^\sqcup$
generated by conjugations and positive Markov moves only.

\proclaim{ Lemma \lemHatBone } {\rm(Key Lemma)}
Let $X\in\hat B_n$, $\eps=\pm1$,
$X'_\eps =X\sigma_n^{-1}\sigma_{n-1}^{2\eps}$ and
$X''_\eps=X\sigma_{n-1}^{2\eps}\sigma_n^{-1}$.
Then $(X'_\eps)_{n+1}\ssim (X''_\eps)_{n+1}$.
\endproclaim

\demo{ Proof } Let $a=\sigma_{n-1}$, $b=\sigma_n$, $c=\sigma_{n+1}$,
$\bar a=a^{-1}$, $\bar b=b^{-1}$, $\bar c=c^{-1}$. Then
\if00{
\def\a{\bar a}\def\b{\bar b}\def\c{\bar c}
\let\u\underline
$$
\split
X'_1=&X \b a\, a \overset{\text{Mm}}\to\longrightarrow 
   X\, \u{\b a b}\, \u{b c\b\b}\, a
  =X   a b\, \u{\a\c\c}\, b\, \u{c a}
  =X   a b \c\c\a b a\, \u c
\\
\overset{\text{cyc}}\to\longrightarrow
  &\u{c X a}\, b \c\c\,\u{\a b a}
  = X a\,\u{c b \c\c}\,   b a\b
  = X a\,\u{\b\b c b b}\, a\b
\overset{\text{Mm}}\to\longrightarrow X''_1\\
X'_{-1}=
   &X\,\u{\b\a\b}\,b\a
  = X     \a\b\,\a   b\a
 \overset{\text{Mm}}\to\longrightarrow 
     X   \a\,\u{\b    \b     c     b}\,\a   b    \a
 =\u{X   \a      c}\,  b    \c    \c   \a   b    \a
 =\u c    X     \a     b    \c    \c   \a   b    \a
\\
\overset{\text{cyc}}\to\longrightarrow
        &X    \a     b\,\u{\c     \c   \a}\,b\,\u{\a    c}
 =        X    \a     b     \a\,\u{\c   \c   b     c}\,\a
 =        X    \a     b     \a\,\u{ b    c  \b}\, \b   \a
\overset{\text{Mm}}\to\longrightarrow
          X    \a     b\,\u{\a     \b    \a}
 =X''_{-1}. \qed
\endsplit
$$}\fi
\enddemo

Let $\deg_q:\hat B_n\to\Z$ be the monoid homomorphism such that $\deg_q(q_i)=1$
and $\deg_q(\sigma_i)=0$
for any $i$. We call $\deg_q(X)$ the {\it $q$-degree\/} of $X$.

\proclaim{ Lemma \lemHatBtwo } {\rm(Diamond Lemma) }
If $(Xq_n)_n\ssim (X'q_m)_m$, then either
$(X\sigma_n^{-1})_{n+1}\ssim (X'\sigma_m^{-1})_{m+1}$
or there exist $Z,Z',Z'',Z'''\in\hat B_\infty^\sqcup$
related to $X\sigma_n^{-1}$ and $X'\sigma_m^{-1}$ as follows
{\rm(}the arrows represent negative Markov $q$-moves
which decrease the $q$-degree{\rm)}:
$$
 \matrix
 X\sigma_n^{-1} & \ssim & Z          &       & Z''' & \ssim & X'\sigma_m^{-1} \\
                &       & \downarrow &       & \!\!\downarrow      \\
                &       & Z'         & \ssim & \!Z''
 \endmatrix
                                                       \eqno(\eqDiagram)
$$
\endproclaim

\midinsert
\epsfxsize=30mm
\centerline{\epsfbox{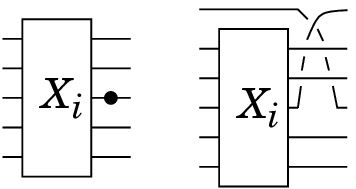}}
\centerline{$Y_i$ \hskip 15mm $Z_i$}
\botcaption{Figure \figYZ}
\endcaption
\endinsert

\demo{ Proof }
Since $Xq_n\ssim X'q_m$, there exists a sequence of words
$Xq_n=Y_0,Y_1,\dots,Y_t$
of the form $Y_i=X_iq_{k_i}\in\hat B_{n_i}$
such that $Y_t$ is a cyclic permutation of $X'q_m$ and
for any pair of consecutive indexes
$i,j$ ($j=i\pm1$) one of the following possibilities holds up to exchange of $i$
and $j$:
\roster
\item"$(i)$" $n_j=n_i$, $k_j=k_i$, $X_i$ and $X_j$ represent the same element of
             $\hat B_\infty^\sqcup$;
\item"$(ii)$" $n_j=n_i+1$, $k_j=k_i$, $X_i=UV$, $X_j=U\sigma_{n_i}V$;
\item"$(iii)$" $n_j=n_i$, $X_i=U\sigma_\ell^\eps$, $X_j=\sigma_\ell^\eps U$,
               $k_j=T_\ell(k_i)$, $\eps=\pm1$;
\item"$(iv)$" $n_j=n_i$, $k_j=k_i\ne\ell$, $X_i=Uq_\ell$, $X_j=q_\ell U$.
\endroster
For $i<j$, we denote $\sigma_i\sigma_{i+1}\dots\sigma_{j-1}$ by $\pi_{i,j}$
and we set $\pi_{i,i}=1$.
Let
$Z_i = X_i\pi_{k_i,n_i}\sigma_{n_i}^{-1}\pi_{k_i,n_i}^{-1}\in\hat B_{n_i+1}$
(see Figure \figYZ). 
It is enough to prove that:
\roster
\item "(a)" $Z_i\ssim Z_j$ in all cases $(i)$--$(iv)$
           (this implies $X\sigma_n^{-1}=Z_0\ssim Z_t$) and
\item "(b)" either $Z_t=X'\sigma_m^{-1}$ or
            we have $Z_t\ssim Z\to Z'\ssim Z''\leftarrow Z'''\ssim X'\sigma_m^{-1}$
            where the arrows mean the same as in (\eqDiagram).
\endroster
Assertion (a) either is evident or follows from Lemma \lemHatBone.
For example,
in Case $(ii)$, we have $Z_i\ssim Z_j$ because 
$$
\split
   Z_i&=UV\pi_{k_i,n_i}\sigma_{n_i}^{-1}\pi_{k_i,n_i}^{-1}
      \ssim U\sigma_{n_i}^{-1}\sigma_{n_i+1}\sigma_{n_i}
      V\pi_{k_i,n_i}\sigma_{n_i}^{-1}\pi_{k_i,n_i}^{-1}
     \overset{\text{\eightpoint def}}\to=
      Z'_i,
\\
   Z_j&=U\sigma_{n_i}V\pi_{k_i,n_i+1}\sigma_{n_i+1}^{-1}\pi_{k_i,n_i+1}^{-1}.
   \qquad\text{and}\quad
   \sigma_{n_j}Z_j\sigma_{n_j}^{-1} = Z'_i.
\endsplit
$$
In Case $(iii)$, $k_i=\ell+1$, $\eps=1$ we have
$Z_i\ssim Z_j$ by Lemma \lemHatBone\ because
$$
\split
    Z_i &= U\sigma_\ell \pi_{\ell+1,n_i} \sigma_{n_i}^{-1} \pi_{\ell+1,n_i}^{-1}
        = V 
          \sigma_{n_i-1}^2\sigma_{n_i}^{-1}
          W, 
\\
    Z_j &= U\pi_{\ell,n_i} \sigma_{n_i}^{-1} \pi_{\ell,n_i}^{-1}\sigma_\ell
         = V
           \sigma_{n_i}^{-1}\sigma_{n_i-1}^2
           W \qquad\text{for}
\\
    V &= U \pi_{\ell+1,n_i}\pi_{\ell,n_i-1}\sigma_{n_i-1}^{-1},
\qquad
    W = \pi_{\ell,n_i-1}^{-1}\pi_{\ell+1,n_i}^{-1}.
\endsplit
$$
In all the other cases Assertion (a) is either similar or easier.

It remains to prove Assertion (b). We know that $Y_t$ is a cyclic permutation of
$X'q_m$. If $Y_t=X'q_m$, then $Z=Z_t=X'\sigma_m^{-1}$ and we are done.
Otherwise we have $X'=U q_k V$ and $Y_t = V q_m U q_k$ for some $k\le m$.
Then we have:
$$
\split
   &Z_t = V q_m U \pi_{k,m} \sigma_m \pi_{k,m}^{-1}
      \ssim
    \sigma_m U \pi_{k,m} \sigma_m \pi_{k,m}^{-1} V \sigma_m^{-1} q_{m+1}
      \overset{\text{\eightpoint def}}\to=
     Z\to Z'
\\
   &X'\sigma_m^{-1} = U q_k V \sigma_m^{-1}
      \ssim
      \pi_{k,m+1}^{-1} V \sigma_m^{-1} U \pi_{k,m+1} q_{m+1}
      \overset{\text{\eightpoint def}}\to=
     Z'''\to Z''.
\endsplit
$$
It is easy to check that $Z'$ and $Z''$ are conjugate. \qed
\enddemo

{\bf Remark \remFT.}
Theorem \thHatB\ admits also a geometric proof based on the interpretation of
the marked points as local modifications introduced in
[\refFT] which increase the Thurston-Bennequin number
(see the extended version of [\refOSh]).


\subhead\S\sectHe. From $A$ to $\hat A$
\endsubhead
Let the notation be as in \S1 and let $\hat A\hat B_\infty$ be the semigroup
algebra of $\hat B_\infty$ with coefficients in $\hat A$.
We have $kB_\infty\subset AB_\infty\subset\hat A\hat B_\infty$.
%
Let $\hat H_\infty$ be the quotient of 
$\hat A\hat B_\infty$ by the bilateral ideal generated by $\ker\pi$ and let
$\hat\pi:\hat A\hat B_\infty\to\hat H_\infty$ be the quotient map.

Let $\hat R$ be the submodule of $\hat A\hat B_\infty$ generated by
all the elements of the form
$$
   XY-YX,\;\;\; X\sigma_n-uX,\;\;\;
     X\sigma_n^{-1} - Xq_n, \;\;\;
     q_{n+1}^a X - v_a X,\;\;\; q_1^a-v_a
$$
with $X,Y\in \hat B_n$ and $n,a\ge1$.
Let $\hat M=\hat H_\infty/\hat\pi(\hat R)$ and let
$\hat t:\hat H_\infty\to\hat M$ be the quotient map.

\proclaim{ Theorem \thMM }
(a). $M$ and $\hat M$ are isomorphic as $A$-modules.
(b). If, moreover, $H_\infty$ is given by $(\eqCubRel)$, then
$\hat M$ is generated by $\hat t(1)$ as an $\hat A$-module.
\endproclaim

\demo{ Proof } (a). Follows from Theorem \thHatB.
(b). Follows from the fact that
$\hat H_{n+1}=\langle q_{n+1}\rangle\hat H_n+\hat H_n\sigma_n\hat H_n
+\hat H_n\sigma_n^{-1}\hat H_n$ where
$\langle q_{n+1}\rangle=\{1,q_{n+1},q_{n+1}^2,\dots\,\}$. \qed
\enddemo

Thus $\hat M=\hat A/\hat I$ where $\hat I$ is the annihilator of $\hat M$.


\subhead\S\sectI. Description of $\hat I$
\endsubhead
In this section we assume that $\hat H_\infty$ is defined by (\eqCubRel).
Let $F_n^+$ (resp.
$\hat F_n$) be the free monoid freely generated by
$x_1^{\pm1},\dots,x_{n-1}^{\pm1}$ (resp. by
$x_1^{\pm1},\dots,x_{n-1}^{\pm1}$, $q_1,\dots,q_n$) and let $\hat A\hat F_n$
be semigroup algebra of $\hat F_n$ over $\hat A$.
We define the {\it basic replacements} as in
[\refOrFu; \S2.1, $(i)$--$(viii)$] and we add to them
\roster
\item"$(ix)$" $x_iq_j\to q_{T_i(j)}x_i$
\endroster
We define $\hat A\Fred_n$ and $\rr:\hat A\hat F_n\to\hat A\Fred_n$
similarly to [\refOrFu; \S2,2] using the replacements $(i)$--$(ix)$.
Then $\hat A\Fred_n$ is the free $\hat A$-module freely generated by the
elements of the form $qX_1X_2\dots X_{n-1}$, $q\in Q_n$,
$X_i\in S_i$ where $S_i$ are as in [\refOrFu; (5)].
We define $\hat\tau_{n}:\hat A\Fred_{n}\to\hat A\Fred_{n-1}$ by setting
$\hat\tau_n(qq_n^a X x_{n-1} Y )=\rr(qXq_{n-1}^a Y)$,
$\hat\tau_n(qq_n^a X x_{n-1}^{-1} Y )=\rr(qXq_{n-1}^{a+1} Y)$,
$\hat\tau_n(qq_n^a X ) = v_a qX$ for $q\in Q_{n-1}$, $X,Y\in F_{n-1}^+$.
We extend $\hat\tau_n$ to $\hat A\hat F_n$ by setting $\hat\tau_n(X)=\hat\tau_n(\rr(X))$
and we define $\hat\tau:\hat A\hat F_\infty\to\hat A\hat F_0=\hat A$ by
$\hat\tau(X)=\hat\tau_1\hat\tau_2\dots\hat\tau_n(X)$ for $X\in\hat A\hat F_n$.

Let $\sh^n$ be the $\hat A$-algebra endomorphism of $\hat A\hat F_\infty$
defined by $\sh\sigma_i=\sigma_{i+n}$, $\sh q_i=q_{i+n}$. We set $\sh=\sh^1$.
For $X\in F_{n+1}^+$, we define $\rho_{n,X}\in\End_{\hat A}(\hat A\Fred_n)$
by setting $\rho_{n,X}(Y) = \hat\tau_{n+1}(X\sh Y)$.

Let $\hat J_4$ be the minimal $\hat A$-submodule of
$\hat A\Fred_4$ which satisfies the conditions
\roster
\item"(J1)"  $\rr(\rr(X_3X_2)X_1)-\rr(X_3\rr(X_2X_1))\in\hat J_4$
for any $X_j\in\sh^{3-j}S_j\setminus\{1\}$, $j=1,2,3$;
\item"(J2)"  $\rho_{4,X}(\hat J_4)\subset\hat J_4$ for any $X\in S_4$.
\endroster

Similarly, let $\hat J_3$ be the minimal $\hat A$-submodule of
$\hat A\Fred_3$ which satisfies
\roster
\item"(J$1'$)"  $q_i\rr(X)
                 -\rr(\rr(X)q_j)\in\hat J_3$ for any
$X=x_2^{\eps_1}x_1^{\eps_2}x_2^{\eps_3}$, $\eps_1,\eps_3\in\{-1,1\}$,
$\eps_2\in\{-1,0,1\}$, $i=1,2,3$, $j=T_2T_1^{\eps_2}T_2(i)$.
\item"(J$2'$)"  $\rho_{3,X}(\hat J_3)\subset\hat J_3$ for any $X\in S_3$.
\endroster

Let $\hat N=\hat A\Fred_2\otimes_{\hat A}\hat A\Fred_2$.
We define $\hat A$-linear mappings $\hat\tau_N:\hat N\to\hat A$ and
$\rho_\delta:\hat N\to\hat N$, $\delta=(\delta_1,\delta_2)\in\{-1,0,1\}^2$,
by setting
$\hat\tau_N(Y_1\otimes Y_2)=\hat\tau(Y_1Y_2)$,
$\rho_\delta(Y_1\otimes Y_2)
= x_1^{\delta_1}\otimes\hat\tau_3((\sh Y_1)x_1^{\delta_1}\sh Y_2)$.
Let $\hat L$ be the minimal $\hat A$-submodule of $\hat N$ satisfying
\roster
\item"(L1)" $\hat\tau_3(x_2^{\eps_1}x_1^{\eps_2}x_2^{\eps_3})\otimes x_1^{\eps_4}
            - x_1^{\eps_2}\otimes\hat\tau_3(x_2^{\eps_3}x_1^{\eps_4}x_2^{\eps_1})
            \in\hat L$ for any $\eps_1,\eps_3\in\{-1,1\}$ and for any
            $\eps_2,\eps_4\in\{-1,0,1\}$;
\item"(L2)"
            $\rho_\delta(\hat L)\subset\hat L$ for any $\delta\in\{-1,0,1\}^2$.
\endroster

\proclaim{ Theorem \thI }
$\hat I = \hat\tau(\hat J_4) + \hat\tau(J_3) + \hat\tau_N(\hat L)$.
\endproclaim

A proof repeats almost word by word the proof of Main Theorem in [\refOrFu]
(we ignore the variables $q_i$ when we define
the weight function on $\hat F_\infty$).

\smallskip

Each of the modules $\hat J_4$, $\hat J_3$,  $\hat N$ is defined
as the limit of an increasing sequence of submodules of a finite rank
$\hat A$-module. Since $\hat A$ is not Noetherian, this does not give yet
a way to compute them. However, we can approximate $\hat A$ by Noetherian
rings $\hat A_d = A[v_1,\dots,v_d]$ and the projections $\pr_d(\hat I)$
can be effectively computed where $\pr_d:\hat A\to\hat A_d$
is the quotient by the ideal $(v_{d+1},v_{d+2},\dots\,)$.
Namely, let $(\hat A\Fred_n)_d$, $(\hat J_4)_d$, $(\hat J_3)_d$,
$(\hat N)_d$, $(\hat L)_d$ be the $\hat A_d$-modules obtained by the above
procedure but with the additional relations $q_i^{d+1}=0$ for any $i$.
Then we have
$\pr_d(\hat I) = \hat\tau(\hat J_4)_d + \hat\tau(J_3)_d + \hat\tau_N(\hat L)_d$
and these modules (at least theoretically) can be computed as limits of increasing sequences
of Noetherian modules.
The rank of $(\hat A\Fred_4)_d$ (where $(\hat J_4)_d$ sits) is equal to
$315(d+1)^4$. We hope that, at least for $d=1$ or $2$, the computations can be
performed in practice.

\smallskip
{\bf Remark \remBetaZero. }
If $\beta=0$ (the
case when the Groebner base of $I_v$ was computed in [\refOrFu]), then
the obtained transversal link invariants a priori cannot detect transversally
non-simple links. Indeed, in this case we have $1=\alpha\sigma_1^{-1}+\sigma_1^{-3}$,
hence $q_1=q_1(\alpha\sigma_1^{-1}+\sigma_1^{-3})=(\alpha\sigma_1^{-1}+\sigma_1^{-3})q_2
=q_2$. Thus $q_1=q_2=q_3=\dots$ whence $v_1=v_2=\dots$ and we obtain
$M=M_v$, $t=t_v$ and $P_t(L) = (v/u)^{(n-e)/2}P_{t_v}(L)$, i.~e., the
invariant $P_t$ reduces to a usual link invariant $P_{t_v}$ and
Thurston-Bennequin number $n-e$.

\smallskip
{\bf Remark \remGB. } By [\refOrGB], all the computations in the huge module
$(\hat A\Fred_4)_d$ can be done in $\Z/m\Z$ for $m$ not very big.


\enddocument